\documentclass[11pt]{article}
\usepackage{amssymb}
\usepackage{latexsym}

\newtheorem{thm}{Th\'eor\`eme}[section]
\newtheorem{prop}[thm]{Proposition}
\newtheorem{lemma}[thm]{Lemme}
\newtheorem{cor}[thm]{Corollaire}
\newtheorem{dfn}[thm]{D\'efinition}

\newtheorem{rmk}[thm]{Remarque}

\newcommand{\reals}{\mathbb R}

\begin{document}

\title{\bf  Les isomorphismes infinit\'esimaux des tissus de codimension 1. } 

\author{ Jean-Paul Dufour  \\}

\date{}
\maketitle

\begin{abstract} We study local ($n+1$)-webs of codimension 1 on a manifold of dimension $n.$ We  give a complete description of their possible Lie algebras of infinitesimal diffeomorphisms. More precisely we show that these Lie algebras are direct products of sub-algebras which are isomorphic to $\frak{sl}(2),$ to the non-commutative 2-dimensional Lie algebra or commutative. We give also a precise limitation of the number of such direct factors  and examples.
\end{abstract}

\section{Introduction.}

Dans ce texte toutes les fonctions que nous introduirons sont suppos\'ees analytiques r\'eelles. Nos r\'esultats sont
 \'egalement vrais dans le cas complexe. La notation $\partial x$ est utilis\'ee pour simplifier ${\partial}\over{\partial x}.$

Tous les r\'esultats de ce travail sont locaux : on travaille au voisinage d'un point de $\reals^n$ que nous appelons l'origine.

 Un $d$-tissu de codimension 1 est une famille de $d$ feuilletages, de codimension 1,  $n$ \`a $n$ transverses. Il est localement donn\'e par ses $d$ int\'egrales premi\`eres $f_1,f_2,\dots,f_d.$ Un tel tissu local sera not\'e $(f_1,f_2,\dots,f_d).$ Il est dit parall\'elisable si  on peut choisir des  coordonn\'ees locales $x_1,\dots ,x_n$  dans lesquelles toutes les int\'egrales $f_i$ sont des fonctions lin\'eaires. Si $d$ est inf\'erieur ou \'egal \`a $n$ le tissu est toujours parall\'elisable ; par contre la  paral\-l\'elisabilit\'e n'est plus  automatique si $d$ est au moins $n+1.$ 
    De fait nous allons travailler ici dans le cas $d=n+1.$ Dans ce cas nous pouvons supposer les coordonn\'ees choisies pour que notre tissu soit
$$W=(x_1,\dots , x_n,f(x_1,\dots , x_n))\ ,$$
o\`u $f$ est une fonction  dont toutes les d\'eriv\'ees partielles $f_{x_i}$ sont non-nulles.
 
Un isomorphisme infinit\'esimal d'un tel tissu est un champ de vecteurs dont le flot pr\'eserve chacun des $n+1$ feuilletages du tissu.

Elie Cartan a \'etudi\'e le cas $n=2,$ voir \cite{C}. Il a montr\'e que l'ensemble des isomorphismes infinit\'esimaux d'un 3-tissu du plan forme une alg\`ebre de Lie de champs de vecteurs de dimension 0, 1 ou 3 ; jamais de dimension 2.  De plus le cas o\`u la dimension est 3 est celui des tissus parall\'elisables.

Les cas o\`u le tissu est parall\'elisable est celui o\`u il est isomorphe \`a $(x_1,\dots ,x_n,x_1+\cdots +x_n).$ Dans ce cas on voit facilement que l'alg\`ebre de ses isomorphismes infinit\'esimaux est de dimension $n+1$ et isomorphe \`a l'ensemble des translations et homoth\'eties infinit\'esimales.

Dans ce travail nous \'etudions le cas des tissus $$W=(x_1,\dots , x_n,f(x_1,\dots , x_n)) ,$$ {\sl non-parall\'elisables}. Nous  d\'eterminons les types d'alg\`ebres de Lie d'iso-morphismes infinit\'esimaux possibles pour de tels tissus. 

 On note $\frak n$ l'alg\`ebre de Lie non-commutative de dimension 2. Notre r\'esultat principal est le suivant.

\begin{thm} \label{th1} L'ensemble des isomorphismes infinit\'esimaux d'un $(n+1)$-tissu de codimension 1 en dimension $n,$ non parall\'elisable, est  une alg\`ebre de Lie de dimension inf\'erieure ou \'egale \`a $n$ et qui est  isomorphe \`a un produit direct de  sous-alg\`ebres isomorphes soit \`a $\frak{sl}(2),$ soit \`a $\frak n,$ soit commutative. Si cette alg\`ebre est commutative, elle est de dimension strictement plus petite que $n.$ On a aussi la contrainte suivante : si on a $S$ facteurs isomorphes \`a $\frak{sl}(2)$ avec $S>1,$ $N$  facteurs isomorphes \`a $\frak n$ et un facteur commutatif de dimension $C$ on a  
$$n\geq 4S+2N+C-1.$$ \end{thm}

\begin{dfn} On note $\frak{g}$  l'alg\`ebre de Lie des isomorphismes infinit\'esimaux du tissu $W.$  \end{dfn}

La preuve  du th\'eor\`eme pr\'ec\'edent utilise  la strat\'egie suivante.

  Les isomorphismes infinit\'esimaux du tissu $W$ sont les champs de vecteurs $\overrightarrow{X}$ tels que
$$\overrightarrow{X}.x_i =\phi_i(x_i),\ \overrightarrow{X}.f(x_1,\dots , x_n)=\phi(f(x_1,\dots ,x_n)),$$
 pour tout $i$ variant de 1 \`a $n$  et o\`u les $\phi_i$ et $\phi$ sont des fonctions d'une variable. Si l'on ne retient que les premi\`eres relations cela impose que $\overrightarrow{X}$ est de la forme
$$\overrightarrow{X}=\sum_{i=1}^nX(x_i)\partial x_i.$$
Autrement dit les composantes $X_i$ sont des fonctions de la seule variable $x_i.$

 Nous choisissons une base de $\frak{g}$ de la forme 
$$(\overrightarrow{X}^1,\dots ,\overrightarrow{X}^m)$$
avec
 $$[\overrightarrow{X}^r,\overrightarrow{X}^s]=\sum_{u=1}^{m}\lambda^{rs}_u\overrightarrow{X}^u$$
o\`u les $\lambda^{rs}_u$ sont des constantes (les {\sl constantes de structure}).

On note 
$$\overrightarrow{X}^j_i=X_i^j(x_i)\partial x_i,$$
Ce champ de vecteurs  ne d\'epend que de la variable $x_i$ et $\overrightarrow{X}^1_i,$ \dots , $\overrightarrow{X}^m_i$ engendrent une alg\`ebre de Lie, de dimension finie, de champs de vecteurs de la droite.

Le point crucial est le fait que la dimension de telles alg\`ebres de Lie est au plus 3. Si c'est 3, elle est isomorphe \`a $\frak{sl}(2)\ ;$ si c'est 2, elle est isomorphe \`a $\frak{n}.$

Le plan de notre travail est le suivant.

Dans la prochaine section nous montrons qu'un tissu de codimension 1 est parall\'elisable si et seulement s'il admet un isomorphisme infinit\'esimal non trivial qui s'annule en un point. On en d\'eduit que, dans le cas non parall\'elisable, la dimension de $\frak g$ est inf\'erieure ou \'egale  \`a $n$. Dans la section suivante nous pr\'ecisons le choix de l'origine. Dans la section qui suit nous montrons que l'on peut choisir les coordonn\'ees locales pour que les composantes des isomorphismes infinit\'esimaux soient des polyn\^omes de degr\'e au plus 2.  

 Ensuite nous travaillons en deux temps. Dans la premi\`ere  \'etape nous d\'ecrivons toutes les alg\`ebres de Lie de dimension inf\'erieure \`a $n$ form\'ees d'isomorphismes infinit\'esimaux du $n$-tissu trivial $(x_1,\dots ,x_n).$  Elles sont isomorphes \`a un produit direct de sous-alg\`ebres isomorphes soit \`a $\frak{sl}(2),$ soit \`a $\frak n,$ soit commutative.

Dans la deuxi\`eme \'etape, nous pr\'ecisons dans quels cas il existe une fonction $f(x_1,\dots ,x_n)$ tels que $(x_1,\dots ,x_n,f(x_1,\dots ,x_n))$ soit un tissu ayant ce type d'alg\`ebres d'isomorphismes infinit\'esimaux.

Cela donnera la preuve de notre th\'eor\`eme.

Dans notre derni\`ere section nous donnons la classification compl\`ete, \`a isomorphisme pr\`es, des 4-tissus de codimension 1 en dimension 3 qui sont non-parall\'elisables et qui admettent une alg\`ebre de Lie d'isomorphismes infinit\'esimaux non-nulle.

Il est probable que la m\'ethode utilis\'ee dans ce texte  pour d\'emontrer notre th\'eor\`eme principal s'\'etende aux cas des $(n+1)$-tissus de codimension $r$ en dimension $nr,$ au moins pour $r=2.$

\section { Conditions pour qu'un tissu de codimension 1  soit parall\'elisable. }

\begin{prop}\label{prop1} 
Si le tissu $W$ admet un isomorphisme infinit\'esimal qui s'annule \`a l'origine, mais qui n'est pas identiquement nul, alors
il est paral\-l\'elisable.\end{prop}

Pour prouver cette proposition
on rappelle la notion de {\bf forme normale} pour un tissu $W$ (voir \cite{DJ}) : apr\`es des  changements de variables ad hoc  la fonction $f$ a la forme normale locale
$$f(x_1,\dots , x_n)=x_1+\cdots  + x_n +\sum_{i<j}x_ix_ja_{ij}(x_1,\dots , x_n)$$
avec $\sum_{i<j}a_{ij}(0,\dots ,0)=0.$ De plus cette forme normale est unique \`a une homoth\'etie pr\`es $(x_1,\dots , x_n)\mapsto ( \lambda x_1,\dots , \lambda x_n).$ 

On voit qu'il n'existe une homoth\'etie non triviale qui laisse $f$ invariante que si la partie non lin\'eaire de la forme normale de $f$ est nulle. S'il existe un isomorphisme infinit\'esimal non trivial qui fixe l'origine, son flot ne peut \^etre form\'e que d'homoth\'eties ; cela prouve que la forme normale de $f$ est $x_1+\cdots  + x_n $ et donc que $W$ est parall\'elisable, achevant ainsi la preuve de la proposition pr\'ec\'edente.

\begin{cor}\label{cor} Si $W$ est non-parall\'elisable $\frak {g},$ son alg\`ebre de Lie d'iso\-morphismes infinit\'esimaux, est au plus de dimension $n.$ \end{cor}
 \noindent Preuve. Raisonnant par l'absurde on suppose que $\frak{g}$ est de dimension $m$ avec $m>n$ et choisissons une de ses  bases $(\overrightarrow{X}^1,\dots ,\overrightarrow{X}^m).$  Les valeurs  $$(\overrightarrow{X}^1(0),\dots ,\overrightarrow{X}^m(0))$$ de ces champs \`a l'origine donnent $m$ champs de vecteurs constants. Comme on a plus de champs  
que de vecteurs, on voit qu'il existe une combinaison lin\'eaire non triviale des $\overrightarrow{X}^j$ qui s'annule \`a l'origine. Le r\'esultat d\'ecoule alors de la proposition pr\'ec\'edente.

Soit $P$ un point de notre domaine de travail. L'orbite de $P$ pour $\frak g$ est l'ensemble des points que l'on peut joindre \`a $P$ en suivant des trajectoires de champs de vecteurs de $\frak g.$ La proposition pr\'ec\'edente a aussi le corollaire suivant.
\begin{cor}\label{corr} Si $W$ est non-parall\'elisable, ces orbites forment un feuilletage (r\'egulier).\end{cor}

\begin{dfn} On note $\cal G$ ce feuilletage.\end{dfn}

\section{Choix de l'origine.}\label{origine}

On consid\`ere un $(n+1)$-tissu $\cal W$, de codimension 1, non-parall\'elisable sur un ouvert de $\reals^n.$ On note ${\cal F}_i,$ pour $i$ variant de 1 \`a $n+1,$ les feuilletages qui le composent. On suppose qu'il admet une alg\`ebre de Lie d'isomorphismes infinit\'esimaux $\frak g$ de dimension $m.$ Etudions le feuilletage associ\'e $\cal G.$

Choisissons un $i$ particulier et des coordonn\'ees locales $x_1,$ ..., $x_n$  pour que le feuilletage  ${\cal F}_i$ ait l'int\'egrale premi\`ere $x_1.$ Si $\overrightarrow{X}=\sum_iX_i\partial x_i$  est dans $\frak g$ alors sa premi\`ere composante $X_1$ n'est fonction que de la variable $x_1.$ On en d\'eduit que l'on a trois cas. Ou bien cette composante est identiquement nulle et les trajectoires de $\overrightarrow X$ sont incluses dans les feuilles de  ${\cal F}_i.$ Ou bien cette composante est partout non-nulle et les trajectoires
 de $\overrightarrow X$ sont transverses \`a toutes les feuilles de ${\cal F}_i.$  Ou bien cette composante est nulle en un point  mais  non-nulle ailleurs et les trajectoires de $\overrightarrow X$ sont transverses \`a toutes les feuilles de ${\cal F}_i$ sauf celle qui contient ce point. 

On en d\'eduit qu'il y a trois types de feuilles de $\cal G$ au voisinage de l'origine. Soit elles sont toutes incluses dans les feuilles de ${\cal F}_i.$  Soit elles sont toutes transverses  aux feuilles de ${\cal F}_i.$ Soit elles sont toutes transverses  \`a toutes les feuilles de ${\cal F}_i$ sauf \`a un certain nombre de feuilles isol\'ees.

Nous choisissons notre origine hors de ces feuilles isol\'ees pour les diff\'erents feuilletages ${\cal F}_i,$ $i$ variant de 1 \`a $n+1.$

\begin{rmk}\label{rem1} Avec un tel choix de notre origine on peut supposer que, pour tout $i,$ ou bien les feuilles de $\cal G$ sont transverses \`a toutes les feuilles de ${\cal F}_i$, ou bien elles sont toutes incluses dans les feuilles  de ce feuilletage.\end{rmk}

Notons que, si $\frak g$ est de dimension $m,$ les feuilles de $\cal G$ ne peuvent \^etre incluses que dans les feuilles de, au plus, $n-m$ feuilletages du tissu.

\section {On peut supposer que les composantes de la base de $\frak g$ sont polynomiales.}

On consid\`ere le tissu $W$ de la section pr\'ec\'edente et l'alg\`ebre de Lie $\frak g$ de ses isomorphismes infinit\'esimaux.

Nous utilisons  les notations de l'introduction.

\begin{dfn} On notera $M$ la matrice de fonctions dont la i-\`eme colonne est form\'ee des $X^j_i(x_i),$ $j$ variant de 1 \`a $m$.\end{dfn}
 
Le corollaire \ref{corr} peut se traduire de la fa\c{c}on suivante.
\begin{lemma} \label{rang} La matrice $M$ est une matrice de rang maximum \`a $n$ colonnes et $m$ lignes avec $m\leq n.$ \end{lemma}

Le r\'esultat essentiel de cette section est le suivant.
\begin{prop}\label{poly}  Quitte \`a faire des changements de variables $x_i\mapsto g_i(x_i)$ ad hoc, on peut supposer que $M$ a tous ses coefficients polynomiaux de degr\'e au plus 2. De plus les colonnes qui ne sont pas identiquement nulles con\-tiennent toutes au moins un coefficient \'egal \`a 1. \end{prop}

Pour prouver cette proposition on commence par remarquer que l'on peut faire des op\'erations de lignes dans $M\ ;$ cela revient \`a modifier la base de $\frak g.$ On peut aussi permuter des colonnes ; cela revient \`a permuter les variables. Par contre on ne peut pas faire des op\'erations de colonnes arbitraires.

Pour deux fonctions $g(t)$ et $h(t)$ de l'unique variable $t$, on note 
$$[g(t),h(t)]=g(t)h'(t)-h(t)g'(t).$$

Cela d\'efinit un crochet de Lie sur les fonctions diff\'erentiables d'une seule variable : ce n'est qu'une autre fa\c{c}on d'\'ecrire le crochet de Lie des deux champs $g(t)\partial t$ et  $h(t)\partial t.$ On notera que deux fonctions $k(t)$ et $l(t)$ commutent pour ce crochet (i.e. $[k(t),l(t)]=0$) si et seulement si l'une est \'egale \`a l'autre multipli\'ee par une constante.

On rappelle les relations
$$[\overrightarrow{X}^r,\overrightarrow{X}^s]=\sum_{u=1}^{m}\lambda^{rs}_u\overrightarrow{X}^u.$$
Elles ont comme corollaire les relations 
$$[\overrightarrow{X}_i^r,\overrightarrow{X}_i^s]=\sum_{u=1}^{m}\lambda^{rs}_u\overrightarrow{X}_i^u$$
pour tout $i$ qui, avec le crochet de fonctions d'une variable d\'efini ci-dessus, deviennent, plus simplement,
$$[X^r_i,X_i^s]=\sum_{u=1}^{m}\lambda^{rs}_uX_i^u.$$

\begin{dfn} On note $\frak{g}_i$ l'alg\`ebre de Lie de champs de la droite engendr\'ee par $\overrightarrow{X}_i^j$, $j$ variant de 1 \`a $m.$\end{dfn}

Par la suite $\frak{g}_i$ pourra \^etre vue soit comme une alg\`ebre de Lie de champs de vecteurs de la droite soit comme une alg\`ebre de Lie de fonctions d'une variable.

 Rappelons le r\'esultat classique suivant.
\begin{lemma} Les sous-alg\`ebres de Lie de dimensions finies de l'espace des champs de vecteurs de la droite sont de dimension 0, 1, 2 ou 3. Apr\`es un bon choix de coordonn\'ee, celles qui contiennent un champ non nul au voisinage de l'origine, sont des sous-alg\`ebres de l'alg\`ebre des champs polynomiaux de degr\'e deux \`a une variable. \end{lemma}

On peut red\'emontrer ce lemme comme suit. On consid\`ere une base
$(g_1(x)\partial x,\dots ,g_N(x)\partial x)$ de notre alg\`ebre. On note $o(g)$ l'ordre de $g(x)\  ;$  c'est dire que le d\'eveloppement de Taylor de $g$ en 0 est de la forme $ax^{o(g)}+bx^{o(g)+1}+\cdots$ avec $a$ non nul. Quitte \`a modifier notre base par des transformations lin\'eaires, on peut supposer 
$$o(g_1)<\cdots <o(g_N).$$ 
Or pour le crochet des fonctions d'une variable d\'efini plus haut,
on la formule \'evidente
$$[x^r,x^s]=(s-r)x^{r+s-1}.$$ Cette formule montre que, pour $N>1,$ $g_{N-1}$ est d'ordre au plus 1. le premier r\'esultat en d\'ecoule. Enfin s'il y a un champ non nul alors $o(g_1)=0$ et un bon changement de variables nous ram\`ene \`a $g_1=1\ ;$  ce qui permet d'achever la preuve de notre lemme.

Il suffit d'appliquer ce dernier lemme pour chaque variable $x_i$ pour achever la preuve de la proposition \ref{poly}.

\section{Cas o\`u $M$ a  des coefficients de degr\'e deux.}\label{deg2}

Dans cette section nous montrons le r\'esultat suivant.

\begin{prop}\label{degre2} On suppose que la matrice $M$ a au moins un coefficient polynomial de degr\'e deux. Quitte \`a faire des op\'erations de lignes et permuter des colonnes, elle a la forme par blocs
$$M:=\left|\matrix{S & 0\cr 0 & A\cr }\right|$$
o\`u $A$ est une matrice \`a coefficients affines  et $S$ est une matrice diagonale par blocs du type
$$
S:=\left| \matrix{ S_1 & 0 &\cdot & \cdot & 0\cr  0 & S_2 & 0 & \cdot & 0 \cr \cdot &\cdot  &\cdot  &\cdot  & \cdot \cr \cdot & \cdot & \cdot & \cdot & \cdot \cr  0 & \cdot &\cdot & 0 &  S_r \cr  }\right|,$$ o\`u chaque bloc $S_i$ est de la forme
$$\left|\matrix{1 & \cdot & \cdot &\cdot & 1 \cr x_{r(i)+1}+c_{r(i)+1}& \cdot &\cdot &\cdot & x_{r(i)+p(i)}+ c_{r(i)+p(i)}  \cr (x_{r(i)+1}+c_{r(i)+1})^2 & \cdot &\cdot &\cdot & (x_{r(i)+p(i)}+ c_{r(i)+p(i)})^2  \cr} \right |$$
o\`u les $c_j$ sont des constantes et $p(i)$ est sup\'erieur ou \'egal \`a 3.
\end{prop}

Sous l'hypoth\`ese de cette proposition, quitte \`a permuter les variables, on peut supposer que $\frak g_1$ est de dimension 3. C'est donc l'alg\`ebre des polyn\^omes de degr\'e 2 en $x_1.$

 En faisant des op\'erations de lignes sur $M$ on peut supposer
$$X^1_1=1,\ X_1^2=x_1,\ X_1^3=(x_1)^2,\ X_1^4=\cdots=X_1^n=0.$$

On en d\'eduit les 3 relations 
$$[\overrightarrow{X}^2,\overrightarrow{X}^1]=-\overrightarrow{X}^1,$$ 
$$[\overrightarrow{X}^1,\overrightarrow{X}^3]=2\overrightarrow{X}^2,$$ 
$$[\overrightarrow{X}^2,\overrightarrow{X}^3]=\overrightarrow{X}^3.$$

\begin{rmk} On retrouve l\`a le r\'esultat classique suivant : l'alg\`ebre des polyn\^omes de degr\'e 2 est isomorphe \`a $\frak{sl}(2)$ ; l'isomorphisme est  celui qui associe au polyn\^ome $a+bx+cx^2$ la matrice
 $$\left( \matrix{
  b/2 & c \cr
 -a&-b/2 \cr
}\right)\ ;$$
de plus $(x^2,x,1)$ est, \`a multiplication par des constantes pr\`es, une base de Chevalley de cette alg\`ebre.
\end{rmk}

Ainsi dans chaque colonne ($i$ variant de 2 \`a $n$) on a 
$$[X^2_i,X^1_i]=-X^1_i,\ [X^2_i,X^3_i]=X^3_i,\ [X^3_i,X^2_i]=2X^2_i.$$

 On a alors deux cas. Ou bien $(X_i^1,X_i^2,X_i^3)$ forme une base ``de Chevalley"  d'une alg\`ebre de Lie isomorphe \`a $\frak{sl}(2)$ (pour le crochet des fonctions d'une seule variable) ou bien $X_i^1,$ $X_i^2$ et $X_i^3$ sont tous les trois nuls, car on a les constantes de structure d'une alg\`ebre de Lie semi-simple.
Supposons que pour un indice  $i>1,$  on soit dans le premier cas.  Un calcul simple montre alors que $X_i^1$ et $X_i^3$ ne peuvent \^etre simultan\'ement nuls \`a l'origine. Si $X_i^1$ est non nul on peut faire un changement de la variable $x_i$ pour avoir $X_i^1=1$ et l'on en d\'eduit l'existence d'une constante $c_i$ telle que l'on ait
$$X^1_i=1,\ X_i^2=x_i+c_i,\ X_i^3=(x_i+c_i)^2.$$

Rappelons maintenant que la matrice $M(0)$ est de rang maximum. Alors ses trois premi\`eres lignes doivent \^etre ind\'ependantes. Cela implique que l'on a au moins 2 colonnes, autres que la premi\`ere, qui sont du type pr\'ec\'edent. Quitte \`a permuter les colonnes nous supposerons que ce sont les trois pre\-mi\`eres. Alors par des op\'erations de ligne on peut se ramener au cas o\`u il existe $p>2$ tel que les $p$ premi\`eres colonnes soient du type pr\'ec\'edent.

Ainsi nous avons montr\'e que l'on peut supposer que la  sous-matrice de $M$ form\'ee par ses trois premi\`eres lignes a la forme
$$\left|\matrix{1 & \cdots & 1 & 0 &\cdots & 0 \cr x_1+c_1& \cdots & x_{p}+ c_{p} & 0 &\cdots & 0  \cr ( x_1+c_1)^2& \cdots & (x_{n}+ c_{p})^2 & 0 &\cdots & 0 \cr } \right |,$$ avec $p\geq 3$ et des constantes $c_i$ telles que $c_1,$  $c_2$ et $c_3$ sont deux \`a deux diff\'erentes.

Montrons maintenant que la sous-matrice de $M$ obtenue en lui enlevant les trois premi\`eres lignes et les $n-p$ derni\`eres colonnes est nulle. 

On sait d\'ej\`a que la premi\`ere colonne de $M$ a tous ses coeficients, autres que les trois premiers, nuls. Etudions la deuxi\`eme  colonne de $M.$ Comme ses coefficients $X^j_2$ doivent satisfaire les m\^emes relations de commutations, vis \`a vis du crochet de Lie, que les $X_1^j,$ on en d\'eduit que les $X_2^j$ pour $j>3$ doivent commuter avec $1,$ $x_2+c_2$ et $ (x_2+c_2)^2$ : on en tire la nullit\'e de tous ces coefficients. On recommence ce raisonnement dans chacune des $p$ premi\`eres colonnes pour prouver la nullit\'e de notre sous-matrice
 
Par les calculs pr\'ec\'edents on s'est ramen\'e au cas 
$$M:=\left|\matrix{S_1 & 0\cr 0 & M'\cr }\right| ,$$
avec 
$$S_1=\left|\matrix{1 & \cdot & \cdot &\cdot & 1 \cr x_{1}+c_{1}& \cdot &\cdot &\cdot & x_{p}+ c_{p}  \cr (x_{1}+c_{1})^2 & \cdot &\cdot &\cdot & (x_p+ c_p)^2  \cr} \right |,$$
o\`u $M'$ est une matrice \`a coefficients polynomiaux de degr\'e au plus deux. 

Si $M'$ est \`a coefficients affines nous avons prouv\'e notre proposition. Si $M'$ a des coefficients de degr\'e deux, on recommence la proc\'edure pr\'ec\'edente en rempla\c{c}ant $M$ par $M'$ : cela fait apparaitre un second bloc $S_2$ et nous ram\`ene au cas d'une matrice $M$ de taille strictement plus petite. Ce processus prouve notre proposition en un nombre fini d'\'etapes.

On a le corollaire \'evident.
\begin{cor} Dans le cas o\`u $M$ a des coefficients de degr\'e deux, l'alg\`ebre $\frak g$ est la somme directe d'un nombre fini (non nul) de sous-alg\`ebres isomorphes \`a $\frak{sl}(2)$ et d'une sous-alg\`ebre $\frak{a}$ correspondant \`a une matrice $A$ \`a coefficients affines. Chaque feuilletage de notre tissu ne peut \^etre transverse qu'aux orbites d'une seule de ces sous-alg\`ebres. \end{cor}

\section{Cas o\`u $M$ est \`a coefficients affines.}

La  proposition \ref{degre2} nous ram\`ene \`a l'\'etude du cas o\`u $M$ n'a que des coefficients affines. Si ces coefficients sont tous constants cela veut dire que $\frak g$ est commutative. C'est le cas le plus simple que nous \'etudierons plus pr\'ecis\'ement plus loin. Dans cette section nous nous int\'eressons d'abord au cas o\`u $M$ a au moins un coefficient non constant.

Dans cette section  nous allons prouver le r\'esultat suivant.
\begin{prop}\label{degre1} On suppose que la matrice $M$ est \`a coefficients affine pas tous constants.  Quitte \`a faire des op\'erations de lignes et permuter des colonnes, elle a la forme par blocs
$$M:=\left|\matrix{N & B\cr 0 & C\cr }\right|$$
o\`u $B$ et $C$ sont \`a coefficients constants  et $N$ est une matrice diagonale par blocs du type
$$
N:=\left| \matrix{ N_1 & 0 &\cdot & \cdot & 0\cr  0 & N_2 & 0 & \cdot & 0 \cr \cdot &\cdot  &\cdot  &\cdot  & \cdot \cr \cdot & \cdot & \cdot & \cdot & \cdot \cr  0 & \cdot &\cdot & 0 &  N_r \cr  }\right|,$$ o\`u chaque bloc $N_i$ est de la forme
$$\left|\matrix{1 & \cdot & \cdot &\cdot & 1 \cr x_{r(i)+1}+c_{r(i)+1}& \cdot &\cdot &\cdot & x_{r(i)+p(i)}+ c_{r(i)+p(i)}  \cr } \right |$$
o\`u les $c_j$ sont des constantes.
\end{prop} 

D\'emontrons cette proposition.

 Quitte \`a permuter des colonnes de $M$ et faire des op\'erations de lignes, on peut supposer que la premi\`ere colonne est 
$$\left|\matrix{1 \cr x_1 \cr 0 \cr \cdot \cr \cdot \cr \cdot \cr 0 \cr } \right |.$$

Cela impose la relation $[\overrightarrow{X}^1,\overrightarrow{X}^2]=\overrightarrow{X}^1$ et donc
$$[X^1_i,X^2_i]=X^1_i$$
pour tout $i.$ On en d\'eduit, utilisant la remarque \ref{rem1}, que toutes les colonnes  sont de l'un des deux types suivants
$$\left|\matrix{1 \cr x_i+c_i \cr 0 \cr \cdot \cr \cdot \cr \cdot \cr 0 \cr } \right |,\ \left|\matrix{0 \cr c_i \cr X_i^3\cdot \cr \cdot \cr\cdot \cr \cdot \cr X_i^m \cr } \right |,$$
o\`u les $c_i$ sont des constantes.

On peut donc supposer que la sous-matrice de $M$ form\'ee par ses deux premi\`eres lignes est
$$N_1=\left|
\matrix{1 & \cdots & 1 & 0 & \cdots & 0 &\cr  x_1+c_1 & \cdots & x_p+c_p & c_{p+1} & \cdots & c_n
 \cr }\right|$$ 
o\`u les $c_j$ sont des constantes. 

On voit ainsi que $M$ est de la forme
$$M:=\left|\matrix{N_1 & B_1\cr 0 & C_1\cr }\right|$$
o\`u $B_1$ est \`a coefficients constants.
Si  $C_1$ n'a que des coefficients constants la proposition \ref{degre1} est prouv\'ee. Sinon des permutations des $n-p$ derni\`eres colonnes et des op\'erations sur les $m-2$ derni\`eres lignes permettent de supposer que la $(p+1)$-i\`eme colonne est de la forme
$$\left|\matrix{0\cr 0\cr1 \cr x_{p+1}+c_{p+1} \cr 0 \cr \cdot  \cr \cdot \cr 0 \cr } \right |.$$
Avec des arguments analogues \`a ceux utilis\'es ci-dessus, on peut supposer que les colonnes de $p+1$ \`a $p+q$ sont toutes du m\^eme type et l'on obtient que $M$ est de la forme

$$\left|\matrix{N_1 & 0 & A_2\cr 0 & N_2 & B_2\cr 0 & 0 & C_2\cr }\right|$$
o\`u $A_2$ et $B_2$ sont des matrices \`a coefficients constants et $C_2$ est \`a coefficient affines avec

$$N_2=\left|
\matrix{1 & \cdots & 1 \cr  x_{p+1}+c_{p+1} & \cdots & x_{p+q}+c_{p+q} 
 \cr }\right|.$$ 

On ach\`eve donc la d\'emonstration de la proposition \ref{degre1} en it\'erant cette m\'ethode.

On en d\'eduit directement le r\'esultat suivant.

\begin{cor} Dans le cas o\`u $M$ est \`a coefficients affines l'alg\`ebre $\frak g$ est somme directe d'un nombre fini de sous-alg\`ebres isomorphes \`a $\frak n$ ou commutative.
\end{cor}

Les sous-alg\`ebres isomorphes \`a $\frak n$ sont engendr\'ees par les deux lignes qui correspondent \`a une sous-matrice $N_i.$ La sous-alg\`ebre commutative est donn\'ee par les derni\`eres lignes. On voit que chacun des feuilletages qui constituent le tissu peuvent \^etre transverses \`a plusieurs des orbites de sous-alg\`ebres isomorphes \`a $\frak n$ (ou commutative).

Les propositions \ref{degre1} et \ref{degre2} montrent que, comme annonc\'e dans le th\'eor\`eme \ref{th1}, $\frak g$ est n\'ecessairement produit direct de sous-alg\`ebres isomorphes \`a $\frak{sl}(2),$ $\frak n$ ou commutatives. Pour achever la preuve de ce th\'eor\`eme il faut montrer que, si l'on se place dans ses hypoth\`eses,  on sait exhiber une $(n+1)$-i\`eme fonction $f=f(x_1,\dots ,x_n)$ telle que $W=(x_1, \dots , x_n,f(x_1,\dots ,x_n))$ ait l'alg\`ebre d'isomorphismes infinit\'esimaux $\frak g$ correspondante.

\section{Exemples avec $\frak g$ commutative.}

Dans les derni\`eres sections de ce texte nous  donnons des exemples de tissus $W= (x_1,\dots ,x_n,f(x_1,\dots ,x_n)), $ non-parall\'elisables qui admettent les types possibles d'alg\`ebres de Lie $\frak g$ d'isomorphismes infinit\'esimaux. On commence par le cas o\`u $\frak g$ est commutative.

\begin{lemma}Il existe des $(n+1)$-tissus $W$ non-parall\'elisables qui poss\`edent des $\frak g$ commutatives de dimension $m$ pour tout $m<n$ mais pas avec $m=n.$ \end{lemma}

 Pour $m<n-1$ peut prendre des exemples du type 
$$f(x_1,\dots , x_n)=x_{1}+\cdots +x_m+g(x_{m+1},\dots x_n)$$
o\`u $g$ est une fonction de $n-m$ variables telle que le tissu $$(x_{m+1},\dots ,x_n,g(x_{m+1},\dots ,x_n))$$ n'ait pas d'isomorphismes infinit\'esimaux.

 Pour obtenir un cas avec $m=n-1$ on peut prendre une fonction $f$ du type
$$x_{1}+\cdots +x_{n-2}+g(x_{n-1}-x_n)$$
avec une fonction $g$ non lin\'eaire.

Si le tissu $W $ admettait $n$ isomorphismes infinit\'esimaux qui commutent, on pourrait suposer que ces isomorphismes soient $\partial x_1, \dots , \partial x_n$ et cela entrainerait  que $f$ serait lin\'eaire.

On remarque que ce dernier lemme est la g\'en\'eralisation du r\'esultat d'Elie Cartan pour la dimension 2.

\section{Exemples avec des facteurs isomorphes \`a $\frak n$}\label{exN}

Rappelons que nous travaillons en dimension $n>2.$

\begin{lemma} Il existe des $(n+1)$-tissus non-parall\'elisables  tels que $\frak g$ est le produit de $N$ facteurs $\frak n$ avec une alg\`ebre commutative de dimension $C$ pour tout $N\geq 1$ et tout $C$ avec $2N+C\leq n.$\end{lemma}

Pour prouver ce lemme on pose
$$\overrightarrow{F}^j=\partial x_{2j}+\partial x_{2j-1},\ \overrightarrow{E}^j=(x_{2j}+1)\partial x_{2j}+x_{2j-1}\partial x_{2j-1},\  \overrightarrow{Z}^k=\partial x_k,$$
on voit que l'on a 
$$[\overrightarrow{F}^j,\overrightarrow{E}^j]=\overrightarrow{F}^j$$
pour tout $j$ variant de 1 \`a $N$ et $k$ variant de $2N+1$ \`a $n,$ et que tous les autres crochets des champs $\overrightarrow{F}^r$, $\overrightarrow{E}^s$ et $\overrightarrow{Z}^k$ sont nuls. Donc l'ensemble de ces champs forme la base d'une alg\`ebre de Lie $\frak h$ de champs de vecteurs isomorphe au produit direct de $N$ sous-alg\`ebres isomorphes \`a $\frak n$ et d'une alg\`ebre commutative de dimension $n-2N.$ On voit aussi que c'est une alg\`ebre de Lie d'isomorphismes infinit\'esimaux du tissu trivial $(x_1,\dots ,x_n)$ et qu'elle est de dimension $n.$

 On consid\`ere maintenant  la fonction $f$ d\'efinie par 
$$f(x_1,\dots ,x_n)=\prod_{j=1}^{N}(1+x_{2j}-x_{2j-1}) \exp( \sum_{s=2N+1}^nx_s)$$
avec $2N$ inf\'erieur ou \'egal \`a $n.$

On a les relations
$$\overrightarrow{F}^j.f=0,\ \overrightarrow{E}^j.f=f,\ \overrightarrow{Z}^k.f=f,$$
pour tous $j$ variant de 1 \`a $N$ et $k$ variant de $2N+1$ \`a $n.$ 

On en d\'eduit que le feuilletage d\'efini par les surfaces de niveau de $f$ est laiss\'e invariant par $\frak h.$

Cela nous donne un exemple de tissu $W$ avec un $\frak g$ isomorphe au produit direct de $N$ copies de $\frak n$ et d'une alg\`ebre commutative de dimension $n-2N.$

On modifie cet exemple facilement pour construire des $(n+1)$-tissus de codimension 1 avec une alg\`ebre d'isomorphismes infinit\'esimaux isomorphe \`a  ${\frak n}^N\oplus {\reals}^S $ avec $N\geq 1$ et $2N+S< n.$ Cela prouve le lemme.

\section{Exemples avec des facteurs isomorphes \`a $\frak{sl}(2).$}

C'est le cas le plus difficile.

 Supposons que $\frak g$ ait $S$ facteurs isomorphes \`a $\frak{sl}(2).$ Cela veut dire qu'elle contient  $3S$ champs de vecteurs
$$\overrightarrow{F}^1,\overrightarrow{H}^1,\overrightarrow{E}^1,\dots ,\overrightarrow{F}^S,\overrightarrow{H}^S,\overrightarrow{E}^S,$$ avec
$$[\overrightarrow{F}^j,\overrightarrow{H}^j]=\overrightarrow{F}^j,\  [\overrightarrow{H}^j,\overrightarrow{E}^j]=\overrightarrow{E}^j,\ [\overrightarrow{F}^j,\overrightarrow{E}^j]=2\overrightarrow{H}^j,$$ 
pour tout $j$ variant de 1 \`a $S,$ les autres crochets deux \`a deux de ces champs \'etant nuls. On doit aussi avoir $3S$ fonctions d'une variable
$$\phi^1,\psi^1,\rho^1,\dots ,\phi^S,\psi^S,\rho^S,$$ telles que
$$\overrightarrow{F}^j.f=\phi^j \circ f,\ \overrightarrow{H}^j.f=\psi^j \circ f,\ \overrightarrow{E}^j.f=\rho^j \circ f.$$  
pour tout $j.$

Toutes ces relations impliquent 

$$[\phi^j ,\psi^j ]=\phi^j ,\ [\psi^j ,\rho^j ]=\rho^j ,\ [\phi^j ,\rho^j ]=2\psi^j,$$
pour tout $j,$ les autres crochets de ces fonctions deux \`a deux \'etant nuls.

On en d\'eduit  qu'il ne peut y avoir que deux cas. Ou bien un seul des triplets $(\phi^j,\psi^j,\rho^j)$ poss\`ede des fonctions non nulles, ou bien toutes ces fonctions sont nulles. Dans le premier cas, et si l'on impose que le triplet non nul est celui d'indice $j=1,$ cela entraine qu'il existe
un diff\'eomorphisme local $\theta$ avec
$$\phi^1\circ\theta (t)=1,\ \psi^1\circ\theta (t)=t+c,\ \rho^1\circ\theta (t)=(t+c)^2$$
avec une constante $c,$ tous les autres $\phi^j,$ $\psi^j$ et $\rho^j$ \'etant nuls. Comme  dans notre tissu $W$ on peut remplacer $f$ par $\theta\circ f,$ cela nous permet de simplifier nos relations. 

Remarquons que  le premier cas est celui o\`u le feuilletage par les surfaces de niveau de $f$ est transverse aux orbites de l'un des facteurs $\frak{sl}(2)$, et alors les autres sont incluses dans des surfaces de niveau de $f.$ Le deuxi\`eme est celui o\`u les orbites des facteurs $\frak{sl}(2)$  sont toutes incluses dans des surfaces de niveau de $f.$

On consid\`ere le facteur $\frak{sl}(2)$ donn\'e par $\overrightarrow{F}^1,\overrightarrow{H}^1,\overrightarrow{E}^1.$ Comme nous l'avons vu dans la section \ref{deg2}, on peut supposer
$$\overrightarrow{F}^1=\sum_{i=1}^p\partial x_i,\ \overrightarrow{H}^1=\sum_{i=1}^p(x_i+c_i)\partial x_i, \  \overrightarrow{E}^1=\sum_{i=1}^p(x_i+c_i)^2\partial x_i,$$ 
avec des constantes $c_i.$ On a d\'ej\`a vu que l'on peut supposer $c_1=0.$ Pour simplifier la r\'edaction des prochains calculs, nous supposerons $c_2=1.$ Il est facile de traiter, de mani\`ere analogue, le cas o\`u $c_2$ est plus g\'en\'eral.

 On d\'efinit  les $p-2$ fonctions $\theta_3,\dots ,\theta_p$ par
$$\theta_j(x_1,x_2,x_j)={{x_j-x_1-c_j(x_2-x_1)}\over{1+x_2-x_1}}$$ 
et on v\'erifie de fa\c{c}on \'el\'ementaire  les deux lemmes suivants.

\begin{lemma}\label{l1} Soit $h(y_3,\dots ,y_p,x_{p+1},\dots ,x_n)$ une fonction de $n-2$ variables qui v\'erifie l'\'equation aux d\'eriv\'ees partielles
\begin{equation}\label{eq1}\sum_{j=3}^p((y_j+c_j)^2-(y_j+c_j))\partial y_jh=h^2-h.\end{equation}
Alors la fonction $f$ d\'efinie en prenant pour $f(x_1,\dots ,x_n)$ la quantit\'e
$$x_1+(1+x_2-x_1)h(\theta_3(x_1,x_2,x_3),\dots ,\theta_p(x_1,x_2,x_p),x_{p+1},\dots ,x_n)$$
v\'erifie les relations 
$$\overrightarrow{F}^1.f=1,\ \overrightarrow{H}^1.f=f,\ \overrightarrow{E}^1.f=f^2.$$\end{lemma}

\begin{lemma}\label{l2} Soit $h(y_3,\dots ,y_p,x_{p+1},\dots ,x_n)$ une fonction de $n-2$ variables qui v\'erifie l'\'equation aux d\'eriv\'ees partielles
\begin{equation}\label{eq2}\sum_{j=3}^p((y_j+c_j)^2-(y_j+c_j))\partial y_jh=0.\end{equation}
Alors la fonction $f$ d\'efinie en prenant pour $f(x_1,\dots ,x_n)$ la quantit\'e
$$h(\theta_3(x_1,x_2,x_3),\dots ,\theta_p(x_1,x_2,x_p),x_{p+1},\dots ,x_n)$$
v\'erifie les relations 
$$\overrightarrow{F}^1.f=0,\ \overrightarrow{H}^1.f=0,\ \overrightarrow{E}^1.f=0.$$\end{lemma}

Via des changements rationnels des variables $y_i$ on peut donner des solutions explicites des \'equations (\ref{eq1}) et (\ref{eq2}). 

Le lemme \ref{l1} montre qu'il existe des exemples de fonctions $f$ pour lesquelles le tissu $W=(x_1,\dots , x_n,f(x_1,\dots ,x_n))$ admet un $\frak g$ contenant un facteur direct isomorphe \`a $\frak{sl}(2)$ de fa\c{c}on que les orbites de ce facteur soient transverses au feuilletage par les surfaces de niveau de $f.$ Il existe m\^eme  des exemples avec $p=3.$  En effet le cas $n=3$ avec $\frak{g}$  isomorphe \`a $\frak{sl}(2)$ a fait l'objet de la publication \cite{D}. Nous rappelons son r\'esultat principal.

\begin{prop}\label{DD} On consid\`ere la fonction
$$f(x,y,z)={{ayz+bzx+cxy}\over{ax+by+cz}}$$
avec $a,$ $b$ et $c$ non nuls et $a+b+c=0.$ Alors $(x,y,z,f(x,y,z))$ est, au voisinage de presque tout point $(x_0,y_0,z_0),$ un 4-tissu qui a une alg\`ebre
 de Lie d'isomorphismes infinit\'esimaux  isomorphe \`a $\frak{sl}(2).$  R\'eciproquement tout 4-tissu de codimension 1 sur $\reals ^3$ qui a une alg\`ebre de Lie d'isomorphismes infinit\'esimaux isomorphe \`a $\frak{sl}(2)$ est localement isomorphe au 4-tissu ci-dessus.\end{prop}

Les calculs sont bas\'es sur le fait que $M$ est une matrice de Van Der Monde.

Le lemme \ref{l2} donne des exemples de fonctions $f$  pour lesquelles le tissu $W$ admet un $\frak g$ contenant un facteur direct isomorphe \`a $\frak{sl}(2)$ de fa\c{c}on que les orbites de ce facteur soient toutes incluses dans les feuilles du feuilletage par les surfaces de niveau de $f.$ Il y a cependant une subtilit\'e dans ce deuxi\`eme cas. En effet, pour avoir une fonction $f$ qui donne un vrai tissu il faut que  les d\'eriv\'ees partielles de $f$ \`a l'origine soient toutes non-nulles. Or si l'on choisit un exemple comme dans le lemme \ref{l2} avec $p=3,$ les trois relations 
$$\overrightarrow{F}^1.f=0,\ \overrightarrow{H}^1.f=0,\ \overrightarrow{E}^1.f=0$$
impliquent que $f$ est ind\'ependante de $x_1,$ $x_2$ et $x_3$. Cela interdit ce cas. Il faut donc imposer $p>3$ pour qu'une $f$ donn\'ee par le lemme \ref{l2} convienne. De tels exemples existent toujours. En particulier
$$f={{3+xy-x-3y+t+3z+zt-zx-yt}\over{4+xy-2x-2y+2t+2z+zt-zy-xt}}$$ 
est une fonction des quatre variables $x,$ $y,$ $z$ et $t$ telle que le tissu $(x,y,z,t,f)$ admette un $\frak g$ isomorphe \`a $\frak{sl}(2)$ avec des orbites incluses dans les surfaces de niveau de $f.$ C'est l'exemple avec $c_1=0,$ $c_2=1,$ $c_3=2$ et $c_4=3.$

Ce dernier exemple a une g\'en\'eralisation du type proposition \ref{DD}.
\begin{prop} On consid\`ere la fonction
 $$g(x,y,z,t)={{xy+zt-xz-yt}\over{xy+zt-zy-xt}}.$$
Alors $(x,y,z,t,g(x,y,z,t))$ est, au voisinage de presque tout point, un 5-tissu qui a une alg\`ebre
 de Lie d'isomorphismes infinit\'esimaux  isomorphe \`a $\frak{sl}(2)$ et dont les orbites sont incluses dans les hypersurfaces de niveau de $g.$
\end{prop}

Cette proposition a une g\'en\'eralisation avec r\'eciproque que nous ne  donnerons pas ici car nous n'en avons pas besoin pour d\'emontrer le th\'eor\`eme \ref{th1}.

Supposons que l'on veuille construire un exemple de tissu $W$ avec un $\frak g$ qui contient plusieurs facteurs isomorphes \`a  $\frak{sl}(2).$ Cela peut se faire en combinant les deux lemmes pr\'ec\'edents (et \`a la condition $n>6$). Pour deux facteurs, un avec orbites transverses, un avec orbites incluses dans les surfaces de niveau de $f,$ on peut prendre pour $f$ une fonction qui associe \`a $(x_1,\dots ,x_n)$ la valeur
$$ x_1+(1+x_2-x_1)h(\theta_3(x_1,x_2,x_3),\theta_6(x_4,x_5,x_6),\theta_7(x_4,x_5,x_7),x_{8},\dots ,x_n) $$
o\`u $\theta_3$ est d\'efini comme juste avant le lemme \ref{l1} mais avec
$$\theta_j(x_4,x_5,x_j)={{x_j-x_4-c_j(x_5-x_4)}\over{1+x_5-x_4}}$$ 
pour $j$  \'egal \`a 6 et 7. On impose alors que la fonction $h$ de variables $y_3,y_6,y_7,x_8,\dots ,x_n,$ v\'erifie les deux \'equations aux d\'eriv\'ees partielles :
$$((y_3+c_3)^2-(y_3+c_3))\partial y_3h=h^2-h,$$
$$\sum_{j=6}^7((y_j+c_j)^2-(y_j+c_j))\partial y_jh=0.$$

On peut rajouter d'autres facteurs $\frak{sl}(2)$ en it\'erant le proc\'ed\'e pr\'ec\'edent. Da fa\c{c}on un peu plus pr\'ecise : si l'on a $n\geq 4S-1$ (i.e. $n\geq 3 +4(S-1)$) on peut construire des exemples de $f$ qui donnent un tissu $W$ dont le $\frak g$ a $S$ facteurs directs  isomorphes \`a $\frak{sl}(2)$.

Enfin si l'on fait le produit d'une telle $f$ avec une des fonctions $f$ de la section \ref{exN}, mais avec des variables diff\'erentes,
on construit des tissus $W$ qui ont tous les $\frak g$ permis par le th\'eor\`eme \ref{th1}. Ceci ach\`eve la preuve de ce th\'eor\`eme.

\section{En dimension 3.}

Dans cette derni\`ere section nous r\'esumons les r\'esultats des sections pr\'ec\'eden\-tes en dimension 3. Cela nous donne une classification, \`a isomorphisme pr\`es, des tissus  non-parall\'elisables $W=(x,y,z,f(x,y,z))$  qui admettent une alg\`e\-bre de Lie $\frak g$ d'isomorphismes infinit\'esimaux non-nulle.

\subsection{Cas o\`u $\frak g$ est commutative.}

{\bf Si $\frak g$ est de dimension 1.} Il y a trois cas : ou bien elle est engendr\'ee par $\partial x,$ ou bien par $\partial x+\partial y,$ ou bien par $\partial x+\partial y+ \partial z.$ Dans le premier cas
 $$f(x,y,z)=x+g(y,z)\ ;$$ 
dans le deuxi\`eme 
 $$ f(x,y,z)=ax+g(y-x,z)\ ;$$
 dans le troisi\`eme
 $$ f(x,y,z)=ax+g(y-x,z-x)\ ;$$
 pour une fonction $g$ g\'en\'erique (il faut imposer des contraintes sur la partie lin\'eaire pour que $W$ soit un vrai tissu et \'eviter qu'il soit parall\'elisable) et $a$ est \'egale \`a 0 ou 1.

{\bf Si $\frak g$ est de dimension 2.} Il y a un seul cas : celui o\`u elle est engendr\'ee par  $\partial x$ et $\partial y+ \partial z.$ On a alors
 $$f(x,y,z)= x+ay+ h(y-z),$$
pour une fonction $h$ g\'en\'erique et $a$ \'egal \`a 0 ou 1.

\subsection{Cas o\`u $\frak g$ est isomorphe \`a $\frak n.$}

On a les diff\'erents sous-cas suivants.

 {\bf Sous-cas 1.} On se place dans le cas o\`u $\frak g$ est engendr\'ee par $\partial x$ et $x\partial x+ \partial y$ et on a 
$$f(x,y,z)= x+\exp (y) h(z),$$
pour une fonction $h$ g\'en\'erique.

 {\bf Sous-cas 2.} On se place dans le cas o\`u $\frak g$ est engendr\'ee par $\partial x$ et $x\partial x+ \partial y+\partial z$ et on a 
$$f(x,y,z)= x+\exp (y) h(z-y),$$
pour une fonction $h$ g\'en\'erique.

 {\bf Sous-cas 3.} On se place dans le cas o\`u $\frak g$ est engendr\'ee par $\partial x+\partial y$ et $x\partial x+(y+b) \partial y$ et on a
$$f(x,y,z)=ax+(b+y-x)h(z),$$
pour une fonction $h$ g\'en\'erique et o\`u $a$ est 0 ou 1.

 {\bf Sous-cas 4.} On se place dans le cas o\`u $\frak g$ est engendr\'ee par $\partial x+\partial y$ et $x\partial x+(y+b) \partial y+\partial z$ et on a
$$f(x,y,z)=ax+\exp (z)h((b+y-x)\exp (-z)),$$
pour une fonction $h$ g\'en\'erique et o\`u $a$ est 0 ou 1.

 {\bf Sous-cas 5.} On se place dans le cas o\`u $\frak g$ est engendr\'ee par $\partial x+\partial y+\partial z$ et $x\partial x+(y+b) \partial y+(z+c)\partial z$ et $f(x,y,z)$ a l'expression
$$ax+(b+(y-x)(1+c)-b(z-x))h((b(z-x)-c(y-x))/(b+(y-x)(1+c)-b(z-x))),$$
pour une fonction $h$ g\'en\'erique et o\`u $a$ est 0 ou 1.

\subsection{Cas o\`u $\frak g$ est isomorphe \`a $\frak n\oplus\reals.$}

On a le seul cas 

$$f(x,y,z)=(1+y-x)\exp (z)$$
avec un $\frak g$ engendr\'e par $\partial x+\partial y,$  $x\partial x+(y+1) \partial y$ et $\partial z.$ 

\subsection{Cas o\`u $\frak g$ est isomorphe \`a $\frak{sl}(2).$}

On a le mod\`ele
$$f(x,y,z)={{(z+c)(y+b)+\lambda (z+c)x-(1+\lambda )x(y+b)}\over{x+\lambda (y+b)-(1+\lambda )(z+c)}}$$
qui se d\'eduit directement de la proposition \ref{DD}. Le param\`etre $\lambda$ est diff\'erent de -1 et 0. Des valeurs diff\'erentes de $\lambda$ donnent des tissus non-isomorphes. Par contre on peut choisir $b$ et $c$ arbitrairement de fa\c{c}on qu'ils soient distincts,  non-nuls ainsi que $\lambda b-(1+\lambda )c.$

\end{document}